\documentclass[12pt]{article}
\usepackage{amssymb}

\makeatletter
 
\@addtoreset{equation}{section}

\makeatother

\renewcommand\thefigure{\thesection.\@arabic\c@figure}

\renewcommand\thetable{\thesection.\@arabic\c@table}

\def\reff#1{(\ref{#1})}
\def\sobre#1#2{\lower 1ex \hbox{ $#1 \atop #2 $ } }
\def\supp{{\rm Supp}\,}

\def\proj{{\rm Proj}\,}

\def\cd{{(d-1)}}
\def\c{{c}}

\begin{document}

\def\E{{\Bbb E}}
\def\P{{\Bbb P}}
\def\R{{\Bbb R}}
\def\Z{{\Bbb Z}}
\def\V{{\Bbb V}}
\def\N{{\Bbb N}}
\def\X{{\cal X}}
\def\Y{{\bf Y}}
\def\G{{\cal G}}
\def\T{{\cal T}}
\def\C{{\C}}
\def\C{{\cal D}}
\def\n{{\bf n}}
\def\m{{\bf m}}
\def\b{{\bf b}}
\def\sqr{\vcenter{
         \hrule height.1mm
         \hbox{\vrule width.1mm height2.2mm\kern2.18mm\vrule width.1mm}
         \hrule height.1mm}}                  % This is a slimmer sqr.
\def\square{\ifmmode\sqr\else{$\sqr$}\fi}
\def\one{{\bf 1}\hskip-.5mm}
\def\liml{\lim_{L\to\infty}}
\def\given{\ \vert \ }
\def\ze{{\zeta}}
\def\be{{\beta}}
\def\la{{\lambda}}
\def\ga{{\gamma}}
\def\th{{\theta}}
\def\proof{\noindent{\bf Proof. }}
\def\rate{{e^{- \beta|\ga|}}}
\def\A{{\bf A}}
\def\a{{\bf a}}
\def\B{{\bf B}}
\def\C{{\bf C}}
\def\D{{\bf D}}
\def\K{{\bf K}}
\def\MM{{\bf m}}

\title{Measures on contour, polymer or animal models.\\  A
probabilistic approach}

\author{Roberto Fern\'andez \\[-3mm] {\it \normalsize Consejo Nacional de
    Investigaciones Cient\'{\i}ficas y T\'ecnicas, Argentina}\\[-3mm] 
{\it \normalsize and Universidade de S\~{a}o Paulo}
\\ 
Pablo A. Ferrari\\[-3mm] {\it \normalsize Universidade de S\~{a}o Paulo} 
\\ 
Nancy L. Garcia \\[-3mm] {\it \normalsize Universidade Estadual de Campinas} 
}

\date{}
\maketitle
\noindent{\bf Abstract. } 
We present a new approach to study measures on ensembles of contours,
polymers or other objects interacting by some sort of exclusion
condition.  For concreteness we develop it here for the case of
Peierls contours.  Unlike existing methods, which are based on
cluster-expansion formalisms and/or complex analysis, our method is
strictly probabilistic and hence can be applied even in the absence of
analyticity properties.  It involves a Harris graphical construction
of a loss network for which the measure of interest is invariant.  The
existence of the process and its mixing properties depend on the
absence of infinite clusters for a dual oriented percolation process
which we dominate by a multitype branching process.  Within the region
of subcriticality of this branching process the approach yields: (i)
exponential convergence to the equilibrium measure, (ii) clustering
and finite-effect properties of the contour measure, (iii) a
particularly strong form of the central limit theorem, and (iv) a
Poisson approximation for the distribution of contours at low
temperature.

\vskip 5truemm


\vskip 2truemm 


\vskip 2truemm 

\noindent {\it Short title:} Measures on contour, polymer or animal models.

\section{Introduction}

Contours were introduced by Peierls (1936) to prove the existence of
a first-order phase transition for the Ising model in 2 or more
dimensions.  His argument, later put on a rigorous mathematical basis
by Dobrushin (1965) and Griffiths (1964), used contours only as an
auxiliary device to estimate spin correlations.  Polymer models, in
the sense of interest here, were introduced later by Gruber and Kunz
(1971). These are abstract general models of possibly extended
objects that interact only by volume exclusion.  They include both
contour ensembles and ensembles formed by the open walks
(``polymers'') or surfaces obtained in high-temperature expansions.
Gruber and Kunz (1971) were the first to treat these models as probability
ensembles of their own, and to ask genuinely probabilistic questions
such as existence and properties of the corresponding probability
measure.  The formalism of cluster 
expansions, whose use in mathematical physics started with a
paper by Glimm, Jaffe and Spencer (1976), soon established itself as
the technique of choice to study these type of systems and questions
[Malyshev (1980), Seiler (1982), Brydges (1984)].  The formalism was
extended by Koteck\'y and Preiss (1986) to objects obeying
generalized exclusion laws defined by \emph{compatibility} relations.
This extension was taken up by Dobrushin (1996, 1996a) who proposed
to call \emph{animal} models to such general systems and introduced 
a new approach to the construction of the expansions.

The cluster-expansion technology, being designed to construct and
study distributions of general systems with exclusions, seems
potentially very useful for probabilists in general.  Nevertheless,
its use has so far remained confined to the mathematical physics
community working in statistical mechanics and quantum field theory. 
This unfortunate situation has been pointed out by Dobrushin (1996a),
who attributed it to two reasons:  (1) ``its analytical and
combinatorial complexity'', and (2) ``the absence \ldots [of a]
systematical exposition oriented to mathematicians''.  He addressed
both issues in his posthumous review, Dobrushin (1996a), where he
presented an exposition geared towards ``probabilistic
interpretations and applications'', based on his new approach that
avoids ``tremendous combinatorial considerations''.  In fact, his
approach does not resort to cluster expansions at all.

In our opinion there is still an additional aspect that explains the
lack of popularity, among probabilists, of this powerful technique:
All the existing formulations transcend the probabilistic framework.
First, the expansions used are a bit unnatural from the
measure-theoretical point of view.  Cluster expansions were, in fact,
originally introduced to control the \emph{pressure} of
gases with exclusions.  The rigorous proof that they converge and
have nice mathematical properties required highly nontrivial
combinatorial estimates which took a reasonable form only with the
insight of Cammarotta (1982) [clearly described in Brydges (1984)
and Pfister (1991)].  The existence of a measure is proven by using
pressure-like expansions for numerator and denominator and cancelling
out terms.  It is, therefore, a rather indirect approach whose
mathematical bottleneck refers to an object ---the pressure--- that
from the probabilistic point of view is just auxiliary.  Dobrushin's
approach, on the other hand, avoids the use of the pressure and the
cluster expansion but its use of complex analysis reveals that the
corresponding hypotheses and results go beyond probability.
Furthermore, for actual computations one needs to go back to the 
traditional approach and its explicit expressions for the correlation
functions (or the pressure). 
 
A second manifestation of ``probabilistic unnaturalness'' is conveyed
by the results themselves.  Indeed, the existing formulations require 
the absolute convergence of the expansions involved.
As a consequence, besides existence and mixing properties, they yield 
\emph{analyticity} of the correlation functions with respect
to different parameters, for instance with respect to the
exponential of minus the inverse temperature.   Though analyticity is
a very nice property to have ---in particular it allows Dobrushin to
produce an amazingly simple proof of the central limit theorem--- it
is also a symptom that these approaches are too strong and not
optimal from the probabilist point of view.  This is not just an
academic remark.  The most interesting recent applications of
cluster-expansion methods fall outside these formulations, as they
involve measures that are known or suspected to have 
non-analytical behavior:  Measures at intermediate temperatures
[Olivieri (1988), Olivieri and Picco(1990)], measures for annealed
disordered systems [von Dreifus, Klein and Perez (1995)], measures
for long-range interactions [Bricmont and Kupiainen (1996)] and
infinite-dimensional Sinai-Ruelle-Bowen measures [Bricmont and
Kupiainen (1997)].  

In this paper we present a novel approach to the study of animal models which
presents a number of advantages regarding these issues.  For concreteness we
discuss the case of usual Peierls contours; a more general treatment will be
presented in Fern\'andez, Ferrari and Garcia~(1998a).  Here is an overview of
the main features of our approach.

\begin{enumerate}
  
\item The approach is purely probabilistic, no cluster expansion or complex
  analysis is involved.  The measure is obtained as the unique stationary
  measure of a Markov process.  The condition of
  validity of our theory is stated in terms of a backwards oriented
  percolation process.  The theory holds when percolation is absent.
  
\item The range of validity of the theory exceeds that of previous approaches
  [see comment after \reff{19a}].  Within this range we obtain all the
  properties yielded by the latter ---existence, uniqueness, exponential
  mixing, central limit theorem--- with one conspicuous, and expected,
  exception: analyticity.  
  
\item We obtain a rather nice version of the central limit theorem [stronger
  than that in Dobrushin (1996a)].  

\item The approach allows us to show that the rescaled distribution of
  contours of a fixed length convergence towards a Poisson process.  We are
  not aware of similar results in the literature.
  
\item The construction constitutes, in fact, a simulation scheme that
  converges to equilibrium exponentially fast.  Hence, it has the potential to
  become a very efficient computational tool.

\item   The avoidance of series expansions for the pressure makes our
  approach more direct to compute general properties of the
  equilibrium measure, but limits its use for the estimation 
  of ``thermodynamic'' quantities.  For instance,  the approach does not 
  seem to be suitable for the study of ``surface corrections'' to the
  presure.  Bounds on these corrections are crucial for several
  applications of contour ensembles [see
  eg.\ Zahradn\'{\i}ck (1984), Borgs and Imbrie (1989)].

\end{enumerate}

In this paper we present a careful statement of these results and a
sketch of their proofs. We aim at providing a streamlined exposition
free of inessential technicalities that may obscure the natural form
of the construction.  Nevertheless, we present enough details for an
educated probabilist to reconstruct most of the missing links. The
full argument will be presented in Fern\'andez, Ferrari and
Garcia~(1998), theretofore referred as FFG.

\section{Contour distribution and loss networks. Results.}

\subsection{Contours}

The contours for the ferromagnetic Ising model with ``$+1$'' boundary
conditions, in dimensions $d\ge 2$, 
are surfaces constructed with
$(d-1)$-dimensional unit cubes ---traditionally known as \emph{plaquettes}---
centered at points of $\Z^{d}$ and perpendicular to the edges of the dual
lattice $\Z^d+({1\over 2},\cdots,{1\over 2})$.  
We shall identify a plaquette
with its center and denote $x\in\gamma$ if the plaquette centered at $x$ is in
$\gamma$.  
Two plaquettes are \emph{adjacent} if they have a common
$(d-2)$-dimensional face.  A collection of plaquettes forms a
\emph{connected} surface if for every two plaquettes $x, y$ one can
find a finite sequence of plaquettes, starting at $x$ and ending at $y$, such
that two consecutive plaquettes of the sequence are adjacent.  A
\emph{closed} surface has every $(d-2)$-dimensional face shared by 2
or 4 plaquettes.  A contour, $\gamma$, is a connected and closed
family of plaquettes.  We say that two contours $\gamma$ and $\theta$ are
\emph{incompatible}, and denote $\gamma\cap\theta\neq\emptyset$, if they have
adjacent plaquettes.  We use the notation $\left|x-y\right|$ for the
minimal number of plaquettes needed to link, in a connected fashion,
$x$ with $y$ (this is also known as ``Manhattan distance'').

For $\Lambda\subset \Z^d$, denote by $\G(\Lambda)$ the
set of contours whose plaquettes have centers in $\Lambda$. A configuration
of contours $\eta\in \N^{\G(\Lambda)}$ is a function that at each contour
$\ga$ assigns a natural number $\eta(\ga)$ indicating the number of contours
$\ga$ present in $\eta$.  The subset $\X(\Lambda) \subset \N^{\G(\Lambda)}$ of
compatible-contour configurations is defined as
\begin{equation}
  \X(\Lambda) = \{ \eta \in \{0,1\}^{\G(\Lambda)}\,;\, 
\eta(\ga) \, \eta(\theta) = 0
  \mbox{ if } \ga \cap \theta \neq \emptyset\} \label{eq:X}
\end{equation}
that is, a configuration of contours is compatible if it contains at most one
copy of each contour and does not contain two intersecting contours. 

For each fixed $\beta\in\R^+$, a parameter usually called the inverse
temperature and for each finite $\Lambda\subset\Z^d$ define the measure
$\mu^\Lambda$ on $\X(\Lambda)$ by
\begin{equation}
  \label{141}
  \mu^\Lambda(\eta) = {\exp\Bigl(-\beta \sum_{\ga:\eta(\ga)=1} |\ga|\Bigr)\over
  Z^\Lambda}
\end{equation}
where $|\ga|$ is the area (=number of plaquettes) of the contour
$\gamma$ and $Z^\Lambda$ is a 
normalization constant making $\mu^\Lambda$ a probability.

\subsection{Loss network of contours}

We introduce a birth-and-death dynamics on the set of compatible
contours.  This process is known in the literature as \emph{loss network}, see
Kelly (1991) and references therein. 

We define the process $\eta^\Lambda_t$ as a Markov process on $\X(\Lambda)$
with generator given by:
\begin{equation}
  A^\Lambda f(\eta) = \sum_{\ga \in \G(\Lambda)} \rate \one \{\eta^{+\ga} \in
  \X(\Lambda)\} [f(\eta^{+\ga}) - f(\eta)] + \sum_{\ga \in \G(\Lambda)}
  \eta(\ga)[f(\eta^{-\ga}) - f(\eta)]
\end{equation}
for $f:\X(\Lambda)\to \R$, where $\one\{\,\cdot\,\}$ denotes the
characteristic function of the set $\{\,\cdot\,\}$ and for $\ga\in
\G(\Lambda)$,
\begin{equation}
\label {2.4}
\eta^{\pm\ga}(\theta) = \left\{ \begin{array}{ll}
                           \eta(\theta) & \mbox{ if $\theta \neq \ga$} \\ 
                           \eta(\ga)\pm1   & \mbox{ if $\theta = \ga$} 
                                  \end{array}
                          \right. 
\end{equation}

% and
% \[ \eta^{-\ga}(\theta) = \left\{ \begin{array}{ll}
%                            \eta(\theta) & \mbox{ if $\theta \neq \ga$} \\ 
%                            \max\{\eta(\ga)-1,0\}   & \mbox{ if $\theta =
%                            \ga$}. 
%                                   \end{array}
%                           \right. \]

It is immediate to check that the measure $\mu^\Lambda$ is reversible for
$\eta^\Lambda_{t}$.  

In terms of loss network language, the above process can be described as
follows.  Consider a network consisting of a finite number of links
represented by plaquettes
with vertices in $\Lambda \subset \Z^{d}$, each link comprising one
circuit. Calls are offered to 
this network along routes $\ga \in \G(\Lambda)$ according to independent Poisson
processes with rate $e^{- \beta|\ga|}$. A call accepted on route $\ga$ holds
all links along this route
%circuits from links (vertices of $\Z^{d}$ connected by edges belonging to
%$\ga$) 
for an exponential holding time with mean 1 and on completion of the service
releases all these circuits simultaneously. All arrival streams and holding
times are mutually independent. A call is accepted along route $\ga \in
\G(\Lambda)$ if $\ga$ is not compatible with other calls already in progress.
Hence $\eta^\Lambda_t = (\eta^\Lambda_t(\ga))_{\ga \in \G(\Lambda)}$ where
$\eta^\Lambda_t(\ga)$ is the number of calls in progress on route $\ga$ at
time $t$, then a call is accepted along route $\ga$ at time $t$ if
\[ \sum_{\ga': \ga' \cap \ga \neq \emptyset} \eta^\Lambda_t(\ga') = 0.\]

We can represent this model as a solution of the following system of
equations: 
\begin{equation}
  \eta^\Lambda_t(\ga) = \eta^\Lambda_0(\ga) + \int_{0}^{t}
  \one\Bigl\{\sum_{\ga': \ga' \cap \ga \neq \emptyset} \eta^\Lambda_{s-}(\ga')
  = 0\Bigr\} d N^+_{\ga}(\rate s) - N^-_{\ga} \Bigl( \int_{0}^{t}
  \eta^\Lambda_s(\ga) ds \Bigr) \label{eq:ngamma}
\end{equation}
where $N^+_{\ga}$ e $N^-_{\ga}$ are independent unit Poisson processes;
$N^+_{\ga}$ creates new contours and $N^-_{\ga}$ destroys them. \\

\subsection{Range of validity of the approach}

Let $\X=\{\eta\in\{0,1\}^{\G(\Z^d)}:\eta(\ga)\eta(\th) = 0 \hbox{ if }
\ga\cap\theta \neq \emptyset\}$. Since $\G(\Z^d)$ is countable, $\X$ is
compact in the product topology. Let $f$ be a continuous function on $\X$.
The infinite-volume loss network on $\X$ has formal generator given by
\begin{equation}
  A f(\eta) = \sum_{\ga \in \G} \rate \one \{\eta^{+\ga} \in \X\}
  [f(\eta^{+\ga}) - 
  f(\eta)] + \sum_{\ga \in \G}\eta(\ga) [f(\eta^{-\ga}) - f(\eta)]
\end{equation}
where $\eta^{\pm\ga}$ was defined in \reff{2.4}. 

We use a graphical construction to show that a sufficient condition for the
existence of a process $\eta_t$ on $\X$ with generator $A$ is
\begin{equation}
    \label{143}
    \la_\be=\sum_{\ga \ni 0} |\ga|\,\rate < \infty.
\end{equation}
Using the fact that $\X$ is compact, abstract nonsense imply that, under
\reff{143}, there exists an invariant measure $\mu$ for $\eta_t$. However, the
way of proving existence is so general that we are not able to show any
further property of this measure.  We remark that, as pointed out by Aizenman,
Bricmont and Lebowitz (1987), \reff{143} defines a ``Peierls'' inverse
temperature,
\begin{equation}
    \label{143.f}
    \be_{\rm P} = \inf\{\be: \la_\be < \infty\}\;,
\end{equation}
above which, with probability one, only a finite number of contours
surround any given site (a fact that, for the Ising model, implies
existence of spontaneous magnetization). The results of this paper, however,
apply to the more limited regime
\begin{equation}
  \label{19}
  \be>\be_M,
\end{equation}
where
\begin{equation}
    \label{19a}
   \be_M=\inf\{\be: \la_\be < 1/\cd\}.
%=\sum_{\ga \ni 0} |\ga| \rate < 1\}
\end{equation}
For the Peierls contours the best estimations of the range of validity of
``traditional'' cluster-expansion approaches follow from
Proposition 5.6 in Dobrushin (1996a), which has been stated in its most
precise form by Lebowitz and Mazel (1997).  As a matter of fact, these authors
present their estimations in a form slightly different to ours:
They consider contours with a given site of the dual lattice in
its interior, rather than contours containing a given plaquette as we do.
The final expressions obtained in these two cases are not directly comparable
because they involve differently-aimed upper bounds.  For a meaningful
comparison we have either to transcribe our approach in terms of interior
sites, or to write theirs in terms of anchoring plaquettes.  The latter policy
leads to a bound 
\begin{equation}
  \label{0f}
  \sum_{\ga\ni 0} e^{\c |\ga|}\,e^{-\beta|\ga|}\leq
  {\c \over d-1}
\end{equation}
for some constant $\c>0$.  This bound can be read off the work of Lebowitz and
Mazel (1997) [who obtain $\c = \be e^{- d\be/4}$], where in fact all the hard
estimates [from their formula (2.7) till the end of their paper] refer to
contours containing a fixed plaquette.  On the other hand, our condition
\reff{19}--\reff{19a} implies
\begin{equation}
  \label{2f}
  \sum_{\ga\ni 0} |\ga|\,e^{-\beta|\ga|} \leq   {1 \over d-1},
\end{equation}
which is strictly weaker than \reff{0f} because $e^x > x$ for $x\ge 1$.

Lebowitz and Mazel show that, defining $\be_{LM}$ as the infimum of $\be$
satisfying \reff{0f},
\begin{equation}
  \label{40f}
  \beta_{LM} \ge 64 {\log d \over d},
\end{equation}
where \reff{2f} plus their counting method, yields
\begin{equation}
  \label{41f}
  \beta_{M} \ge 6 {\log d \over d}.
\end{equation}
On the other hand, Aizenmann, Bricmont and Lebowitz (1987) show that the
Peierls temperature defined by \reff{143.f} satisfies
\begin{equation}
  \label{41f}
  \beta_{P} \ge {\log d \over 2d}.
\end{equation}
These three temperatures mark, therefore, limits where different
properties can be proven by perturbation arguments. For $\be\ge
\be_P$, each site of the dual lattice is sorrounded by a finite number
of contours. In spin language, this means lack of percolation of
minority spins (which, in turn, implies symmetry breaking and, by FKG,
non zero magnetization). For $\be\ge \be_M$, in addition, properties
${\bf R1}$---${\bf R5}$ listed below can be proven by
cluster-expansion-like methods. Finally when $\be \ge \be_{LM}$
methods of this type also yield analytic temperature dependence.

\subsection{Results}
We say that $f$ has \emph{support} in $\Upsilon\subset \Z^d$ if $f$ depends
only on contours intersecting $\Upsilon$ (not necessarily contained in
$\Upsilon$). Let $|\supp(f)|= \min\{|\Upsilon|: f$ has support in
$\Upsilon\}$. When we write $\supp(f)$ we mean any $\Upsilon$ such that
$|\Upsilon|=|\supp(f)|$ and $f$ has support in $\Upsilon$. For instance, if
$f(\eta) = \eta(\ga)$, $\supp(f)$ may be set as $\{x\}$ for any $x\in\ga$.

A closer analysis of the graphical construction allows us to show
that for $\be>\be_M$ the following results hold. These are our main
results.

\begin{itemize}
\item[{\bf R1.}] Reversibility and uniqueness: there exists a unique invariant measure $\mu$
  for $\eta_t$.  Furthermore, $\mu$ is reversible for the process $\eta_t$.
  
  \item [{\bf R2.}] The rate of convergence to the invariant measure is exponential. Let
    $\delta_\xi S(t)$ be the distribution of the process at time $t$ when the
    initial configuration is $\xi$. For measurable $f$ we prove
    \begin{equation}
      \label{334}
      \sup_{\xi\in\X} \vert \mu f - \delta_\xi S(t) f\vert \le
      \|f\|_\infty\,|\supp(f)|\, e^{-M_0t}
    \end{equation}
    for any $M_0<(1-\cd\lambda_\be)/(2-\cd\lambda_\be)$.
    
  \item [{\bf R3.}] Infinite-volume limit: Let $\Lambda$ be a (finite or
    infinite) subset of $\Z^d$ and $f$ a measurable function depending on
    contours contained in $\Lambda$.  Then
    \begin{equation}
      \label{130}
      |\mu f - \mu^\Lambda f| \le \|f\|_\infty\,M_2\, \sum_{x\in
       \supp(f)}  e^{-M_3 \, d(x,\Lambda^c)} 
%  |\mu f - \mu^\Lambda f| \le \|f\|_\infty\,(M_2)^2\, \sum_{\scriptstyle
%   x\in \supp(f),\atop\scriptstyle y\in \Lambda^c} |x-y| e^{-M_3 |x-y|}
    \end{equation}
    where $M_2 = e^{(\be - \be_M)/\cd}$ and $M_3  \ge ( \be - \be_M
    )/(d-1)$.  We denoted   $d(x,\Lambda^c) = \min \{|x-y|:
    y\in\Lambda^c\}$.

%    $$
%    M_3  =\left[( \be - \be_M )/\cd + \sum_{\ga \ni 0} e^{\be_M
%        |\ga|/\cd}(1 - e^{(\be - \be_M)|\ga|/\cd})\right]d^{-1/2}.
%    $$
%    
  \item [{\bf R4.}] Clustering. For measurable functions $f$ and $g$ depending
    on contours contained in an arbitrary set $\Lambda\subset\Z^d$:
    \begin{eqnarray}
      \label{101}
      \vert \mu^\Lambda (f g) - \mu^\Lambda f\, \mu^\Lambda g \vert \le
      2\,\|f\|_\infty\,\|g\|_\infty\,(M_2)^2 \,\sum_{\scriptstyle x\in
        \supp(f),\atop\scriptstyle y\in \supp(g)} |x-y| \, e^{-M_3
      |x-y|}\\[-12pt] 
      \nonumber
    \end{eqnarray}
    where $M_2$ ad $M_3$ are the same of \reff{130}. This includes the
    infinite-volume measure $\mu^{\Z^d}=\mu$.
  
  \item [{\bf R5.}] Central limit theorem. Let $f$ be a measurable function on
    $\X$ with finite support such that $\mu f = 0$ and $\mu(
    |f|^{2+\delta})<\infty$ for some $\delta>0$. Assume $D= \sum_x
    \mu(f\tau_xf)>0$. Then $D<\infty$ and
    \begin{equation}
      \label{90}
      {1\over \sqrt {|\Lambda|}}\, \sum_{x\in \Lambda} \tau_x f
      \sobre{\textstyle\Longrightarrow}{\textstyle\Lambda\to\Z^d} \hbox{\rm
      Normal} (0,D)
    \end{equation}
    where the double arrow means convergence in distribution. This result
    generalizes (the central limit) Theorem 7.4 of Dobrushin (1996a). In the
    latter, only functions depending on a finite number of contours are
    considered. 
\end{itemize}

For the following result we write $\mu_\be$ to stress the $\beta$ dependence
of $\mu$.

\begin{itemize} 
\item [{\bf R6.}] Poisson approximation. Let $\eta^\beta$ distributed with
  $\mu_\beta$.  For each measurable $V\subset \R^d$ let
\begin{equation}
  \label{121}
  V(a) = \{x\in \Z^d: x/a\in V\}.
\end{equation}
For each $j$, the process $N^{j,\beta}$ defined by
\begin{equation}
\nonumber
N^{j,\beta}(V) = \sum_{\ga \subset V(e^{\be j}), |\ga|=j} \eta^\beta(\ga).
\end{equation}
converges weakly to a unit Poisson process on $\R^d$ as $\beta \rightarrow
\infty$. The rate of convergence is exponential in $\beta$.

\end{itemize}

The key to the proof of the above results is a graphical construction of the
process starting from a marked Poisson process in $\Z^d\times\R$.  The marks
determine random cylinders whose bases are the contours and the heights are
exponentially distributed random times. The exclusion condition is imposed
through the study of the ``ancestors'' of each cylinder (Section \ref{S3}).
These ancestors determine a (backwards) oriented percolation process, and our
construction is feasible if there is no such percolation.  This is the meaning
of the condition $\be>\be_M$.  All our results follow from the estimation of
the spatial and temporal extension of the cluster of a (finite number of)
cylinders(s) (Section \ref{S4} and \ref{s131}).  This estimation is done
through a domination of the percolation process by a multiple branching
process which, in the regime $\be>\be_M$ has exponential moments (Section
\ref{S6}). The proof of {\bf R6}, also based in the above properties, is
omitted here. A complete proof is presented in FFG.

Ferrari and Garcia (1998) used space-time percolation to show ergodicity of
loss networks under low arrival-rate of calls.

\section {Graphical construction. The BO-cluster}
\label{S3}

\subsection{Finite volume}
To each contour $\ga \in \G$ we associate a Poisson process of rate $\rate$,
and to each time event $T_k(\gamma)$ of the Poisson process we associate an
independent exponentially distributed time $S_k(\gamma)$ of mean one. The
collection $\C=(T_k(\ga), S_k(\ga))_{\ga\in\G,k\in\Z}$ is a family of
double-sided independent marked Poisson processes, with the convention
$T_{-1}(\ga)<0<T_0(\ga)$. The $k$th attempt of birth of a contour $\ga$
occurs at time $T_k(\ga)$; $S_k(\ga)$ corresponds to the lifetime of the
contour. Each triplet $(\ga,T_k(\ga), S_k(\ga))$ is called a cylinder of basis
$\ga$ birth-time $T_k(\ga)$ and lifetime $S_k(\ga)$. To each contour $\th$
present in the initial configuration $\eta_0=\eta$ we independently associate
an exponential time $S(\th)$ and cylinder $(\th,0,S(\th))$. The collection of
initial cylinders is called $\C(0)$. We realize the dynamics $\eta_t$ as a
(deterministic) function of $\C$ and $\C(0)$. 

When the number of possible contours is finite, the construction for $t>0$ is
as follows.  We construct inductively $\K_{[0,t]}$, the set of \emph{kept}
cylinders. The complementary set corresponds to \emph{erased} cylinders. First
include all cylinders of $\C(0)$ in $\K_{[0,t]}$. Then, move forward in time
and consider the first Poisson mark: The corresponding cylinder is erased if
it intersects any of the cylinders already in $\K_{[0,t]}$, otherwise it is
kept.  This procedure is successively performed mark by mark until all
cylinders born before $t$ are considered. Define $\eta_t\in\X(\Lambda)$ as
\begin{eqnarray}
  \label{12}
 \eta_t(\gamma)  &=& \eta_0(\ga)\, \one\{S(\ga)>t\} +\one\{ \exists k\,:\,
  (T_k(\ga),T_k(\ga)+S_k(\ga))\ni t\hbox{ and } (\ga,T_k(\ga), S_k(\ga))
  \hbox{ is kept}\},\nonumber \\
\end{eqnarray}
that is, $\eta_t$ signals all contours which are basis of a kept
cylinder that is alive at time $t$.

It is tedious but easy to show that $\eta_t$ has as generator an operator
defined as $A$, but with the sums restricted to the finite set of contours
involved. In particular, when the contours are contained in a finite
region $\Lambda$, we obtain the process $\eta^\Lambda_t$ with generator
$A^\Lambda$.

The above finite-volume construction can also be performed in
$(-\infty,\infty)$.  Indeed, $\eta^\Lambda_t$ is an irreducible Markov process
in a \emph{finite} state space. Hence, with probability one there exists a
sequence of ordered random times $t_k(\C)$ such that no cylinder in $\C$ is
alive by time $t_k(\C)$. Furthermore $\E(t_{k+1}-t_k)<\infty$.  Therefore one
can apply the above construction independently in each interval
$[t_k(\C),t_{k+1}(\C))$. In this case the cylinders of $\C(0)$ play no role.
This procedure is time-translation invariant and so is the distribution of
$\eta_t$.  This distribution is precisely given by the measure $\mu^\Lambda$.

\subsection{Infinite volume}
For infinite volume, the Poisson processes are indexed by an infinite
set of contours. Hence, it is not possible to decide which is the first mark
in time. The construction must be performed more carefully. There are two
alternatives.

The first alternative is to divide the time interval $[0,t]$ in successive
intervals of small length $h$ and perform the construction in each one of
those intervals.  Under \reff{143} and for small $h$, it is possible to
partition $\Z^d$ in finite regions such that each contour born in $[0,h]$ is
contained in exactly one of these regions.  To show this, one considers the
percolation of (projected) contours and dominates the area occupied by the
contours by a branching process. Such a construction is at the heart of Harris
(1972) original graphical construction of particle systems and it is reviewed
by Durrett (1995).  The mark-by-mark construction described above can be
performed in each of these finite regions to construct the process in the time
interval $[0,h]$. The same procedure can then be applied in the interval
$[h,2h]$, etc.

The second alternative is the one we really use. In order to know whether a
cylinder $C\in\C$ is kept, one has to look at the set of cylinders $C'$ (born
before $C$ and) alive at the birth-time of $C$ whose basis intersects the
basis of $C$. This set is called the first generation of \emph{ancestors} of
$C$. The second generation of ancestors of $C$ consists, previsibly, of the
ancestors of the ancestors, that is those cylinders that are in the first
generation of ancestors of some $C'$ in the first generation of ancestors of
$C$.  Recursively we construct in this way the $n$th generation of ancestors
of $C$.  The set of ancestors (of any generation) of $C$ is called the
\emph{BO- (backwards oriented) cluster} of $C$ and it is denoted by $\A(C)$.
This set may contain cylinders in $\C(0)$.  [We remark that this BO-cluster is
a cluster of \emph{space-time cylinders}, it is different from the usual
cluster of contours considered in the classical works on cluster expansions.]
If for some $C$ the BO-cluster of $C$ has a finite number of cylinders, then
we can decide whether $C$ is kept or not by looking at $\A_{[0,t]}(C) =
\{C'\in\A(C): C'$ is born in $[0,t]\}$, the set cylinders in the BO-cluster of
$C$ born in $[0,t]$.  This is done in the following way.  First, those
$C'\in\A_{[0,t]}(C)$ that have no ancestors are kept.  Then we look to the
remaining cylinders in $\A_{[0,t]}(C)$ and erase those that have a kept
cylinder in its first generation.  We repeat these two steps for the cylinders
in $\A_{[0,t]}(C)$ that have not already declared to be kept or erased, and
continue in this way until we reach $C$.  The end result is a partition of
$\A_{[0,t]}(C)$ in two subsets formed, respectively, by kept and erased
cylinders.  In particular, the subset to which $C$ belongs decides its status.
In fact, under \reff{143} one can prove that all cylinders have a finite
number of ancestors born in the interval $[0,t]$, and, thus, the process
$\eta_t$ can be constructed following, BO-cluster by BO-cluster, the steps of
the finite case \reff{12}.  In the next section, we sketch the proof of the
finiteness of the number of ancestors and give further details of the
construction.

It is natural, and it turns out to be convenient, to
extend the notion of ancestors of a cylinder to that of ancestors of a
space-time point $(x,t)$, $x\in\Z^d$ and $t\in\R$:  Let the first
generation of ancestors of $(x,t)$ be the set of cylinders in $\C$
whose basis contains $x$ and are alive at time $t$.  The $n$th
generation of ancestors of $(x,t)$ is then formed by the $(n-1)$-th
generation of ancestors of the cylinders in the first generation.
The union of all the generations of ancestors is the BO-cluster
$\A(x,t)$ of $(x,t)$.  More generally, the set of ancestors 
of $\Upsilon\subset \Z^d$ at time $t$ is defined by
\begin{equation}
  \label{167}
  \A(\Upsilon,t)=\bigcup_{x\in\Upsilon}\A(x,t).
\end{equation}

\section{Existence of $\mu$ and exponential convergence}
\label{S4}

\subsection{Backwards percolation}
To perform the construction described in the previous section, every cylinder
$C\in\C$ must have a finite number of ancestors. If this is the case, we say
that there is no \emph{(backwards oriented) percolation} in $\C$. Hence, if
with probability one there is no backwards percolation, the double infinite
construction holds and we have a process $(\eta_t)_{t\in\R}$ that is
time-translation invariant. The marginal distribution of $\eta_t$ does not
depend on $t$ and it is called $\mu$. By construction $\mu$ is an invariant
measure for $\eta_t$. This shows the existence of $\mu$ in a constructive way.
In contrast, the existence of $\mu$ under \reff{143} uses a fixed point
theorem.

The condition $\be>\be_M$ implies that there is no percolation with
probability one. This is shown by dominating the number of plaquettes in the
bases of the cylinders in a BO-cluster by a branching process. The number
\begin{equation}
  \label{30f}
  (d-1)\,\la_\be
\end{equation}
is an upperbound on the mean number of branches of the process. That is, the
mean number of plaquettes born from the branching (= incompatible contours) of
each single plaquette. The process is subcritical if this number is less than
one, thus the condition $\be>\be_M$.  This argument, sketched in Section
\ref{S6} below is inspired by Hall (1985), who dominated a continuum
percolation process by a branching process. We sketch this domination in
Section \ref{S6} below. If there is no percolation, the number $|\A(x,t)|$ is
finite for all $(x,t)$. As a consequence, there exists a function $\Phi
:(f,\A(\supp(f),t))\mapsto \Phi(f,\A(\supp(f),t)$ such that for any $f$ with
finite support
\begin{equation}
  \label{77}
  f(\eta_t) = \Phi(f,\A(\supp(f),t)).
\end{equation}
For instance, to decide whether a contour $\ga$ is present at time
$t$ it suffices to look at the BO-cluster of $(x,t)$ for some
$x\in\ga$. The function $\Phi$ is 
the one that decides which cylinders are kept and indicates the
presence/absence of $\ga$ at time $t$. 

\subsection{Time length and space width of the BO-cluster}
Most of the stated properties ---uniqueness of $\mu$, exponential clustering
and finite-volume effects, and exponential convergence to equilibrium of the
loss network--- follow from the observation that for $\be>\be_M$ both the time
length and the space width of the BO-cluster of any given site decay
exponentially. More precisely, let us introduce $\proj(\A(x,t))\subset\Z^d$,
the spatial projection of the BO-cluster, defined as the set of points in
$\Z^d$ belonging to the basis of some cylinder in the BO-cluster:
\begin{equation}
  \label{172}
  \proj(\A(x,t))= \bigcup_{\ga\in\a(x,t)}\{x\in\ga\}
\end{equation}
where $\a(x,t)= \{\ga:(\ga,T_k(\ga),S_k(\ga))\in \A(x,t)$ for some
$k\}$, is the set of bases of the cylinders of the BO-cluster of $(x,t)$.  The
cardinality of this set will be bounded by the cumulative number of points:
\begin{equation}
  \label{171}
  \|\A(x,t)\|=\sum_{\ga\in\a(x,t)} |\ga|\;.
\end{equation}
Indeed, it is clear that 
\begin{equation}
|\proj(\A(x,t))|\le \|\A(x,t)\|\;.
\label{172a}
\end{equation}

We then have:
\begin{enumerate}

\item Let $E_2(t)$ be the set of $\C$ for which the BO-cluster of
$(0,0)$ has time-length larger than $t$:
\begin{equation}
  \label{37}
  E_2(t) = \left\{\C: C\hbox{ is alive at time }-t \hbox{ for some
      }C\in\A(0,0)\right\}.
\end{equation}
Then, for $\be>\be_M$
\begin{equation}
  \label{43}
  \P(E_2(t))\le M_1 \,e^{-t(1-\cd\la_\be)} 
\end{equation}
with $M_1>0$.

\item
Let $E_3(w)$ be the set of $\C$ for which the projection of the 
BO-cluster of $(0,0)$ is not contained in $[-w, w]^d$:
\begin{equation}
  \label{38}
  E_3(w) = \Bigl\{\C:\proj(\A(0,0))\not\subset[-w, w]^d \Bigr\}.
\end{equation}
Then, for $\be>\be_M$ 
\begin{equation}
  \label{44}
  \P(E_3(w))\le M_2 e^{-M_3 w} 
\end{equation}
where $M_2$, $M_3$ are as in \reff{130}.
\end{enumerate}

The proof of \reff{43} and \reff{44} are sketched in Section \ref{S6}.  To
prove \reff{43} we dominate $\A(0,0)$ by a continuous-time branching
process.  On the other hand, to prove \reff{44} we dominate $\|\A(0,0)\|$ by
the total population of a branching process.

\subsection{Proof of R1 and R2}

The exponentially fast time-convergence \reff{334} is a consequence
of \reff{43} and \reff{44}. We use the same Poisson marks to construct
simultaneously the stationary process $\eta_t$ and a process starting at time
zero with an arbitrary initial configuration $\xi$. The second process is
called $\xi_t$, where $\xi_0=\xi$.  The process $\xi_t$ ignores the cylinders
in $\C$ with birth-times less than $0$ and considers cylinders in $\C(0)$ with
basis in $\xi$ and birth-time zero. The process $\eta_t$ ignores the
cylinders in $\C(0)$. Hence for any $\ga\ni 0$,
\begin{equation}
  \label{173}
  |\eta_t(\ga) - \xi_t(\ga)| \le \one \{\A(0,t) \ne
  {\widetilde\A}(0,t)\}  
\end{equation}
where ${\widetilde\A}(0,t)$ is the cluster constructed in $\C_{[0,t]}\cup
\C(0)$.  In FFG it is shown, via a coupling
argument, how \reff{38} and \reff{44} imply that the expectation of
the right hand side of \reff{173} 
decays as $\exp(-M_0\, t)$. The exponential decay of length and width of the
cluster implies exponential decay of the probability that it contains a
cylinder of the initial configuration $\C(0)$.

The uniqueness of $\mu$ follows immediately from \reff{334}.

Reversibility follows from the facts that $(\eta^\Lambda_t)$ converges in
distribution to $(\eta_t)$, $\mu^\Lambda$ converges to $\mu$ and
$\mu^\Lambda$ is reversible for $\eta^\Lambda_t$. From the construction, under
$\be>\be_M$, it is possible to show that $(\eta^\Lambda_t)$ converges almost
surely to $(\eta_t)$. Some details are given in the next sections.

\section{Space-time mixing and the central limit theorem}\label{s131}

\subsection{The key facts} 

The mixing properties of the measure $\mu$ are a consequence of the following
space-time mixing properties of $\C$.

\begin{itemize}
\item Let $f$ be a function depending on contours contained in a finite
  set $\Lambda$. Let $\eta^\Lambda_t$ be the loss network process constructed
  in $\Lambda$. Then
\begin{eqnarray}
  {|\E(f(\eta_0))- \E f(\eta^\Lambda_0)|}\,\le\,
 2\, \|f\|_\infty\,\P\Bigl(\A(\supp(f),0)\neq\A^\Lambda(\supp(f),0)\Bigr)
  \,.\label{73a} 
\end{eqnarray}
where $\A^{\Lambda}(\supp(f),t)$ is the cluster of $(\supp(f),t)$
constructed with cylinders in
\begin{equation}
  \label{133}
  \C^\Lambda = \{(\ga,T_k(\ga),S_k(\ga)) \in\C: \ga \subset \Lambda, k\in\Z\},
\end{equation}
the subset of cylinders whose basis is in $\Lambda$.  

\item For arbitrary measurable functions $f$ and $g$,
\begin{eqnarray}
  \lefteqn{|\E(f(\eta_0)g(\eta_0))- \E f(\eta_0)\,\E g(\eta_0)|}\nonumber\\
  &\le& 2\,\|f\|_\infty\|g\|_\infty \P\Bigl(C'\cap C''\ne \emptyset\hbox{ for
    some } C'\in\A(\supp(f),0) \hbox{ and }C''\in{\widehat\A}
  (\supp(g),0)\Bigr)\nonumber\\ \ \label{73b}
\end{eqnarray}
where ${\widehat\A}(\supp(g),t)$ has the same distribution as $\A(\supp(g),t)$
but is independent of $\A(\supp(f),t)$.
\end{itemize}

The proof of \reff{73a} follows rather straightforwardly
from the space-time construction.   Using \reff{77} we get
\begin{eqnarray}
  f(\eta_0)-  f(\eta^\Lambda_0) &=& 
  \Bigl[\Phi(f,\A(\supp(f),0))-
  \Phi(f,\A^{\Lambda}(\supp(f),0))\Bigr]\nonumber\\
&&\qquad {}\times \one \bigl\{\A(\supp(f),0)\neq\A^\Lambda(\supp(f),0)\bigr\}\;.
\label{81a}
\end{eqnarray}
As, by definition,
$|\Phi(f,\A(\supp(f),t))| \le \|f\|_\infty$,
%The set
%$\Bigl\{\Bigl[\A(\supp(f),0)\Bigr]\,
%\cap\,\Bigl[\A(\Lambda^c,0)\Bigr] = \emptyset\Bigr\}$ is contained
%in 
%$\{\A(\supp(f),0)=\A^\Lambda(\supp(f),0)\}$.
taking expectations and absolute values in \reff{81a} we get \reff{73a}.

The proof of \reff{73b} is similar in
spirit but requires a somewhat more delicate argument based on the coupling of
two continuous-time versions of the backwards percolation process. See details
in FFG.

\subsection{Proof of R3 and R4}

To prove the finite-volume effects \reff{130} we use the space-time
representation \reff{12} and get
\begin{equation}
  \label{63a}
  \mu f - \mu^\Lambda f
=\E f(\eta_0) - \E f(\eta^\Lambda_0).
\end{equation}
By \reff{73a} it is enough to bound 
\begin{equation}
  \label{136}
  \P\Bigl(\A(\supp(f),0)\neq\A^\Lambda(\supp(f),0)
\Bigr),
\end{equation}
which as in \reff{44} is bounded by
\begin{equation}
  \label{165}
  M_2\, \sum_{x\in\supp(f)} e^{-M_3 \, d(x,\Lambda^c)} \;.
\end{equation}
This proves the decay stated in \reff{130}. 

The proof of exponential mixing \reff{101} is similar but using instead the
bound \reff{73b}. 

While we have not yet done a careful study, we believe that \reff{130} and
\reff{101} lead  to sharper inequalities than those obtained via the use of
``duplicated variables'' [von Dreifus, Klein and Perez (1995), Bricmont and
Kupiainen (1996)]. The reason is that clusters formed by superposition of two
systems of contours have larger probabilities of intersection than our
single-system clusters.

\subsection{Proof of the central limit theorem}

We use the results for stationary mixing random fields of Bolthausen (1982).
Let $X_x = \tau_x f$. By hypothesis, $||X_x||_{2+\delta}<\infty$.  Under this
conditions, Bolthausen (1982) shows that if
\begin{equation}
  \label{93}
  \sum_{n=1}^\infty n^{d-1} (\alpha_{2,\infty}(n))^{\delta/(2+\delta)} <\infty
\end{equation}
then $D<\infty$ and \reff{90} holds. Here $\alpha_{2,\infty}(n)$ measures the
dependence between functions depending on the sigma algebra generated by $X_0$
and $X_y$ and the sigma algebra generated by $\{X_x:x\in\Lambda\}$ for
$|\Lambda|=\infty$ and $\min\{|x|, |y-x|: x\in \Lambda\}>n$. In FFG we use
\reff{101} to show that 
\begin{eqnarray}
  \label{120}
  \alpha_{2,\infty}(n) &\le& \,(M_2)^2\,|\supp(f)|\,\sum_{|y|\ge n-2|\supp(f)|} e^{-M_3 |y|}
\end{eqnarray}
Hence, $\alpha_{2,\infty}(n)$ decreases exponentially fast with $n$. This
shows the central limit theorem.

\bigskip

\section{Length and width of the BO-cluster}\label{S6}

To conclude, let us sketch the arguments behind the bounds \reff{43} and
\reff{44}. In both cases we rely on dominating branching processes. 

\subsection{Time length}
To show \reff{43} we consider a continuous time multitype Markov branching
process $\b_t$ on $\N^\G$. 
In this process, each contour $\ga$ lives an mean-one exponential
time after which it dies and gives birth to $k_\th$ contours $\th$,
$\th\in \G$, with probability
\begin{equation}
  \label{30}
  \prod_\th {e^{\mu(\ga,\th)} \left(\mu(\ga,\th)\right)^{k_\th} \over k_\th!}
\end{equation}
for $k_\th\ge 0$. These are independent Poisson distributions of mean
$\mu(\ga,\th)=\one\{\gamma\cap\theta\neq\emptyset\}\,e^{-\beta|\theta|}$. 

Fix $\b_0(\ga) = \left|\{k: (\ga,T_k(\ga),S_k(\ga))\hbox{ is alive at time
    }0\}\right| \,\one\{\ga \ni 0\}$ and zero otherwise. Under this initial
condition it is possible to couple $(\b_t)_{t\ge 0}$ and $\A(0,0)$ in such a
way that
\begin{equation}
  \label{31}
  E_2(t)\, \subset\,\Bigl\{\sum_\th \b_t(\th) = 0\Bigr\}.
\end{equation}
Using the backwards Kolmogorov equation for $R_t = \E\sum_\th\b_t(\th)$, one
can show that
  \begin{equation}
    \label{34}
   \P\Bigl(\sum_\th \b_t(\th) >0\Bigr) \le  e^{(\cd\la_\be-1)t}.  
  \end{equation}

\subsection{Space width}

Define a Galton-Watson branching process $Z_n\in\N$ as follows.
Let $Y^n_i$ be i.i.d.\ non negative integer valued random variables with the
same distribution as
\begin{equation}
  \label{49}
  Y:=\sum_{\gamma \ni 0} \cd|\gamma| X_\gamma
\end{equation}
where $X_\gamma$ are independent integer valued random variables with Poisson
distribution of mean $e^{-\beta |\gamma|}$.
Define $Z_0 = 1$ and 
\begin{equation}
  \label{50}
  Z_{n+1} = \sum_{i=1}^{Z_n} Y^n_i
\end{equation}
(with the convention $\sum_{i=1}^{0} Y^n_i=0$).
It is possible to couple the BO-cluster $\A(x,t)$ and $(Z_n)_{n\ge 0}$ in such
a way that the number of plaquettes in the bases of the cylinders in the $n$th
generation of ancestors of $(x,t)$ is less than or equal to $Z_n$:
\begin{equation}
  \label{169}
   \|\A(x,t)\|\le \sum_{n\ge 0} Z_n.
\end{equation}
Hence, to show \reff{44} it suffices to prove
\begin{equation}
  \label{60}
  \P(Z>k) \le M_2 e^{-M_3 k}
\end{equation}
where $Z=\sum_{n\ge 0}Z_n$. Call $F(b)$ the generating function of $Z$, we
will prove that if $\be>\be_M$, $\bar b = \sup \{b: F(b)<\infty\}>1$.

The generating function of 
$Y$ is given by
\begin{equation}
  \label{51}
  f(a) = \E a^{Y} = \prod_{\ga\ni 0} \E a^{\cd|\ga|X_\ga} =
\exp \left(\sum_{\gamma \ni 0}
  e^{-\be|\ga|}(a^{\cd|\ga|}-1)\right). 
\end{equation}
The radius of convergence of $f(a)$ is given by $\exp(\beta-\beta_{P})$, where
$\beta_P$ is defined in \reff{143.f}.
For $\be>\be_M(>\be_P)$, the radius of convergence is strictly
larger than 1. The mean number of offsprings $\E Z_1$ is given by
\begin{equation}
  \label{52}
  \left. f'(a)\right|_{a=1} = \E Y = \cd\,\sum_{\gamma \ni 0} |\ga|
  e^{-\be|\ga|}<1 
\end{equation}
for $\be>\be_M$. Hence, our branching process is subcritical,
\begin{equation}
  \label{54}
  1 = f(1) \hbox{  and } x = f(x) \hbox{ implies } x\ge 1
\end{equation}
i.e. the smallest solution of the equation $x=f(x)$ is $1$.

By (13.3) of Harris (1963) $F(b)$, the generating function of
$Z$, must satisfy the equation
\begin{equation}
  \label{53}
  F(b) = b f(F(b)).
\end{equation}
The largest solution of this is 
\begin{equation}
  \label{57}
  \bar b= \bar a/ f(\bar a)
\end{equation}
where $\bar a$ is the solution of 
\begin{equation}
  \label{55}
  f'(a) = {f(a) \over a}.
\end{equation}

In this case, it is easy to see that
\begin{equation}
  \label{55.1}
  f'(a) = {f(a) \over a} \cd\,\sum_{\gamma \ni 0} |\ga|  e^{-\be|\ga|} 
a^{\cd|\ga|}
\end{equation}
and $\bar a$ is the solution of
\begin{equation}
  \label{55.2}
\sum_{\gamma \ni 0} |\ga|  e^{-\be|\ga|} a^{\cd|\ga|} = 1/\cd
\end{equation}
which gives us 
\begin{equation}
  \label{55.3}
  \bar{a} = e^{(\be - \be_{M})/\cd}
\end{equation}
Therefore, 
\begin{equation}
  \label{56.2}
  \bar{b} = \exp \Bigl\{ {\be - \be_M\over \cd} 
+ \sum_{\ga \ni 0} e^{\be_M |\ga|/\cd}(1 -
  e^{(\be - \be_M)|\ga|/\cd})\Bigr\}.
\end{equation} 

By exponential Chebichev, fixing $M_2 = \E \bar b^Z $ and $M_3 = \log
\bar b$, we get \reff{60}. 

\bigskip
\baselineskip 12pt
\parskip 1mm  
\section*{Acknowledgments} 
We thank Aernout C. D. van Enter, Joel L. Lebowitz, Fabio Martinelli, Enzo
Olivieri and Roberto Schonmann for some nice discussions.  

This work was partially supported by FAPESP 95/0790-1 (Projeto Tem\'atico
``Fen\^omenos cr\'\i ticos e processos evolutivos e sistemas em equil\'\i
brio'') CNPq, FINEP (N\'ucleo de Excel\^encia ``Fen\'omenos cr\'\i ticos em
probabilidade e processos estoc\'asticos'' PRONEX-177/96).

\vskip 6mm
\baselineskip 13pt
\parskip 0pt
\obeylines
Roberto Fern\'andez, \hfill Pablo A. Ferrari
IEA USP, \hfill IME USP, 
Av. Prof. Luciano Gualberto, \hfill Caixa Postal 66281, 
Travessa J, 374 T\'erreo \hfill 05389-970 - S\~{a}o Paulo,
05508-900 - S\~{a}o Paulo,\hfill BRAZIL
BRAZIL \hfill email: {\tt pablo@ime.usp.br}
email: {\tt rf@ime.usp.br} \hfill http://www.ime.usp.br/\~{}pablo
\vskip 5mm

Nancy L. Garcia
IMECC, UNICAMP, Caixa Postal 6065, 
13081-970 - Campinas SP 
BRAZIL
email: {\tt nancy@ime.unicamp.br}
http://www.ime.unicamp.br/\~{}nancy

\end{document}